\documentclass[20pt,leqno]{article}
\usepackage[left=2.5cm,right=2.5cm]{geometry}
\usepackage{amsmath,amsthm,amscd,amssymb}
\usepackage{latexsym,graphics,longtable,lscape}
\usepackage{epstopdf}
\usepackage{float}
\renewcommand{\baselinestretch}{1}
\newtheorem{remark}{Remark}[section]
\newcommand{\bremark}{\begin{remark} \em}
\newcommand{\eremark}{\end{remark} }

\begin{document}
\parindent 15pt
\renewcommand{\theequation}{\thesection.\arabic{equation}}
\renewcommand{\baselinestretch}{1.15}
\renewcommand{\arraystretch}{1.1}
\def\disp{\displaystyle}
\title{\bf\large Minkowski problem arising from sub-linear elliptic equations
\footnotetext{ \\ 2010 Mathematics Subject Classification: 31B25,35B45,35J96,35J61,52A20.\\
E-mail addresses:  qiuyidai@aliyun.com, 895224263@qq.com\\
\quad The work was supported by NNSFC (Grant: No.11671128).}
\author{{\small Qiuyi Dai, Xing Yi}\\
{\small LCSM, College of Mathematics and Statistics, Hunan Normal University,}\\
 {\small Changsha, 410081, Hunan, People's Republic of China.}\\
 }
}

\date{}

\maketitle

\abstract
{\small Let $\Omega\subset\mathbb{R}^N$ be a bounded convex domain with boundary $\partial\Omega$ and $\nu(x)$ be the unit outer vector normal to $\partial\Omega$ at $x$. Let $\mathbb{S}^{N-1}$ be the unit sphere in $\mathbb{R}^N$. Then, the Gauss mapping $g: \partial\Omega\rightarrow \mathbb{S}^{N-1}$, defined almost everywhere with respect to surface measure $\sigma$, is given by $g(x)=\nu(x)$. For $0<\beta<1$, it is well known that the following problem of sub-linear elliptic equation
\begin{equation*}
\left\{
\begin{array}{ll}
-\Delta\varphi=\varphi^{\beta}&x\in \Omega\\
\varphi>0&x\in\Omega\\
\varphi=0& x\in\partial\Omega.
\end{array}
\right.
\end{equation*}
has a unique solution. Moreover, it is easy to prove that each component of $\nabla\varphi(x)$ is well-defined almost everywhere on $\partial\Omega$ with respect to $\sigma$. Therefore, we can assign a measure $\mu_{\Omega}$ on $\mathbb{S}^{N-1}$ such that $d\mu_{\Omega}=g_*(|\nabla\varphi|^2d\sigma)$. That is
$$\int_{\mathbb{S}^{N-1}}f(\xi)d\mu_{\Omega}(\xi)=\int_{\partial\Omega}f(\nu(x))|\nabla\varphi(x)|^2d\sigma$$
for every $f(\xi)\in C(\mathbb{S}^{N-1})$. The so-called Minkowski problem associated with $\mu_{\Omega}$ ask to find bounded convex domain $\Omega$ so that $\mu_{\Omega}=\mu$ for a given Borel measure $\mu$ on $\mathbb{S}^{N-1}$. Our main results of this paper are the weak continuity of $\mu_{\Omega}$ with respect to Hausdorff metric and the uniquely solvability of Minkowski problem associated with $\mu_{\Omega}$.

\vskip 0.2in
{\bf Key words:} Sub-linear elliptic equation; Minkowski problem; Weak continuity, Borel measure.

\section*{1. Introduction}

\setcounter{section}{1}

\setcounter{equation}{0}

\noindent

Let $N\geq 2$. A non-empty compact, convex subset of $\mathbb{R}^N$ with non-empty interior is called a convex body. Let $\mathcal{K}^N$ be the set of all convex bodies of $\mathbb{R}^N$ and $\mathcal{K}^N_0$ be the set of all convex bodies of $\mathbb{R}^N$ containing the origin in their interior. An open connected subset of $\mathbb{R}^N$ is called a domain. Therefore, the closure $\overline{\Omega}$ of a bounded convex domain $\Omega$ is in $\mathcal{K}^N$ and conversely the interior $K_{\text{int}}$ of a convex body is a bounded domain. Let $\mathbb{S}^{N-1}$ be the unit sphere in $\mathbb{R}^N$. For a domain $\Omega$, we denote by $\partial\Omega$ the boundary of $\Omega$, by $\nu(x)$ the outer unit normal at $x\in\partial\Omega$, by $g: \partial\Omega\rightarrow \mathbb{S}^{N-1}$ the Gauss mapping defined by $g(x)=\nu(x)$. Let $\Omega$ be a bounded domain in $\mathbb{R}^N$, we consider
\begin{equation}\label{eq101}
\left\{\begin{array}{ll}
-\Delta\varphi=\psi(x,\varphi) & x\in \Omega\\
\varphi>0 &x\in\Omega\\
\varphi=0 & x\in\partial\Omega.
\end{array}
\right.
\end{equation}
Let $\varphi(x)$ be a solution to problem (\ref{eq101}). If the gradient $\nabla\varphi$ of $\varphi$ is well-defined almost everywhere on $\partial\Omega$ with respect to surface measure $\sigma$ and $\nabla\varphi\in L^\infty(\partial\Omega)$ , then for $p\geq 0$ one can define a measure $\mu_{p,\Omega}$ on $\mathbb{S}^{N-1}$ so that
$d\mu_{p,\Omega}=g_*(|\nabla\varphi|^p d\sigma)$, or in other words, for any $f(\xi)\in C(\mathbb{S}^{N-1})$, there holds
$$\int_{\mathbb{S}^{N-1}}f(\xi)d\mu_{\Omega}(\xi)=\int_{\partial\Omega}f(\nu(x))|\nabla\varphi(x)|^2d\sigma.$$
The so-called Minkowski problem associated  with $\mu_{p,\Omega}$ can be stated as follows,
\vskip 0.1in
{\em{\bf Problem (M).} For a given non-negative finite Borel measure $\mu$ on $\mathbb{S}^{N-1}$, under what conditions does there exist a bounded convex domain $\Omega$ of $\mathbb{R}^N$ such that
$\mu_{p,\Omega}=\mu$?}
\vskip 0.1in
In the case $p=0$, $d\mu_{0,\Omega}=g_*(d\sigma)$. Therefore, problem (M) is reduced to the classical Minkowski problem which is a central problem in the theory of convex bodies and has been studied extensively in the past years. See for example Minkowski \cite{M1897, M1903}, Aleksandrov \cite{A1996}, Fenchel-Jessen \cite{FJ1938}, Busemann \cite{B1958}, Bonnesen-Fenchel \cite{BF1987}, Lewy \cite{L19382}, Nirenberg \cite{N1953}, Cheng-Yau \cite{CY}, Pogorelov \cite{P78}, Caffarelli \cite{Caff1989,Caff1990sc,Caff1990i,Caff1991s} and so on.

In the case $p=1$ and $\psi(x,s)=\delta(x)$ a Dirac measure at origin, problem (M) is reduced to the prescribing harmonic measure problem which was posed and solved by D. Jerison \cite{J1991}.

In the case $p=2$ and $\psi(x,s)=2$, problem (M) is reduced to the Minkowski problem induced by torsional rigidity which was studied by Colesanti-Fimiani \cite{CF}.

In the case $p=2$ and $\psi(x,s)=\lambda s$, problem (M) is reduced to the Minkowski problem induced by the principle eigenvalue of Dirichlet Laplacian which was studied by D. Jerison \cite{J1996adv}.

In the present paper, we focus our attention on the case $p=2$ and $\psi(x,s)=s^\beta$ with $\beta\in(0,1)$. More precisely, we prescribe a finite Borel measure $\mu$ on $\mathbb{S}^{N-1}$ and find a bonded domain $\Omega$ with $\overline{\Omega}\in\mathcal{K}^N$ such that the measure $\mu_{2;\Omega}$ induced by the unique solution to the following problem
\begin{equation}\label{eq102}
\left\{\begin{array}{ll}
-\Delta\varphi=\varphi^\beta & x\in \Omega\\
\varphi>0 &x\in\Omega\\
\varphi=0 & x\in\partial\Omega.
\end{array}
\right.
\end{equation}
verifies $\mu_{2;\Omega}=\mu$. In the sequel, we always denote $\mu_{2,\Omega}$ by $\mu_{\Omega}$ for simplicity.
\vskip 0.1in
So far, there are two frames used to attack problem (M). One is from the view point of differential equations, the other is from the view point of variational method. The frame from the view point of differential equations is to finding solution of problem (M) firstly in smooth case. This is equivalent to finding $0<h(\xi)\in C^{2,\eta}(\mathbb{S}^{N-1})$ and $0<\varphi(x)\in C^{2,\eta}(\Omega_h)$ for a given smooth positive function $f(\xi)$ on $\mathbb{S}^{N-1}$ such that
\begin{equation}\label{eq104}
\left\{
\begin{array}{ll}
|\nabla\varphi(g^{-1}(\xi))|^p\det(\partial_{ij}h(\xi)+\delta_{ij}h(\xi))=f(\xi) &\xi\in \mathbb{S}^{N-1}\\
-\Delta\varphi=\psi(x,\varphi)&x\in \Omega_h\\
\varphi>0&x\in\Omega_h\\
\varphi=0& x\in\partial\Omega_h\\
\end{array}
\right.
\end{equation}
with $\Omega_h=\{x\in\mathbb{R}^n: x\cdot\xi<h(\xi)\ \ \text{for any}\ \ \xi\in\mathbb{S}^{N-1}\}$. Then, by approximating general non-negative Borel measure by Borel measures with positive smooth density function, the solution of problem (M) for general non-negative Borel measures can be obtained as a limit of solution obtained in the smooth case. The variational frame requires to finding a functional $F(\Omega)$ of domains so that the measure $\mu_{p,\Omega}$ can be generated by the domain variation of $F(\Omega)$. When this is done, the solution of problem (M) can be obtained as a critical point of $F(\Omega)$. No matter what frame is used, the most important tool for solving problem (M) is the weak continuity of the measure $\mu_{p,\Omega}$ with respect to the Hausdorff metric which should be established for any concrete problem. 

To state our results precisely, we recall firstly some notations and definitions of convex bodies. Let $K,L\in\mathcal{K}^N$ and denote by $d(x,y)$ the Eucleadian distance between points $x$ and $y$. The Hausdorff distance $d_H(K,L)$ between $K$ and $L$ is defined by
$$d_H(K,L)=\max\{\sup\limits_{x\in K}\inf\limits_{y\in L}d(x,y), \sup\limits_{x\in L}\inf\limits_{y\in K}d(x,y)\}.$$
Hausdorff distance between to convex bodies can be also expressed by the so-called support function of convex body. For a convex body $K\in\mathcal{K}^N$, the support function $h_K$ of $K$ is defined as
$$h_K: \mathbb{S}^{N-1}\to\mathbb{R},\quad h_K(\xi)=\sup\limits_{x\in K}x\cdot\xi.$$
Geometrically, the quantity $h_K(\xi)$ is the signed distance from the origin of the supporting hyperplane to $K$ having $\xi$ as outer normal. In particular, if $K\in\mathcal{K}_0^N$, then we always have $h_K(\xi)\geq 0$ for any $\xi\in\mathbb{S}^{N-1}$. Let $B_r(x)=\{y\in\mathbb{R}^N: d(y,x)<r\}$. If the diameter $R_\Omega$ of a domain $\Omega$ is defined by
$$R_\Omega=\sup\limits_{x,y\in\Omega}d(x,y),$$
Then, for any $K\in\mathcal{K}^N$, we always have have
$$|h_K(\xi)|\leq R_K \quad\forall \xi\in\mathbb{S}^{N-1}.$$
A more important result we need is the following relationship between the Hausdorff metric and the support function which can be found in \cite{S2014}.
\begin{equation}\label{eq103}
\begin{array}{ll}
d_H(K,L)=\sup\limits_{\xi\in\mathbb{S}^{N-1}}|h_k(\xi)-h_L(\xi)|& \forall K, L\in\mathcal{K}.
\end{array}
\end{equation}
One of our main results is the following weak continuity of $\mu_{\Omega}$.
\vskip 0.1in
{\em {\bf Theorem 1.1:} Assume that $\Omega$ is a bounded convex domain and $\{\Omega_n\}_{n=1}^\infty$ is a sequence of bounded convex domains. If $d_H(\overline{\Omega}_n,\overline{\Omega})\to 0$ as $n\to\infty$, then $\mu_{\Omega_n}$ converges weakly to $\mu_{\Omega}$ in the sense of measure. That is
$$\lim\limits_{n\to\infty}\int_{\mathbb{S}^{N-1}}f(\xi)d\mu_{\Omega_n}=\int_{\mathbb{S}^{N-1}}f(\xi)d\mu_{\Omega}\ \ \ \forall\ \ f\in C(\mathbb{S}^{N-1}).$$}

By making use of Theorem 1.1 and the variational method, we can prove the solvability of Minkowski problem (M) arising from problem (\ref{eq102}). This gives the following Theorem.
\vskip 0.1in
{\em {\bf Theorem 1.2:} Let $\mu$ be a positive finite Borel measure on $\mathbb{S}^{N-1}$ such that
$$\int_{\mathbb{S}^{N-1}}|\theta\cdot\xi|d\mu(\xi)>0\quad\forall\theta\in\mathbb{S}^{N-1}.$$
Then there is a bounded convex domain $\Omega$ such that $\mu_{\Omega}=\mu$ if and only if $\int_{\mathbb{S}^{N-1}}\xi d\mu(\xi)=0$. Moreover, $\Omega$ is uniquely determined up to a translation.}
\vskip 0.05in
Let $\varphi_\Omega(x)$ be the unique solution to problem (\ref{eq102}). Define the domain functional $F(\Omega)$ by the formula
$$F(\Omega)=\int_\Omega|\nabla\varphi_\Omega(x)|^2dx.$$
If $\Omega$ is bounded convex domain of class $C^2$, a result of Y. Huang, C. Song and L. Xu given in \cite{HSX2018} implies that $\mu_\Omega$ can be induced by the variation of $F(\Omega)$. More precisely, if $K$ and $L$ are $C^2$ bounded convex domains and $K_t=K+tL$ represent the Minkowski sum of $K$ and $tL$, then one has
$$\frac{d}{dt}F(K_t)|_{t=0}=\frac{1+\beta}{1-\beta}\int_{\mathbb{S}^{N-1}}h_L(\xi)d\mu_K.$$
With the aid of Theorem 1.1, we can prove that the above fomula is true for general bounded convex domains. For a given positive finite Borel measure $\mu$, we define the mean width $W_\mu(\Omega)$ of $\Omega$ associated with $\mu$ by the formula
$$W_\mu(\Omega)=\frac{2}{\mu(\mathbb{S}^{N-1})}\int_{\mathbb{S}^{N-1}}h_\Omega(\xi)d\mu.$$
Noting that $W_\mu$ is linear under the Minkowski addition, we have
$$\frac{d}{dt}W_\mu(K_t)|_{t=0}=\frac{2}{\mu(\mathbb{S}^{N-1})}\int_{\mathbb{S}^{N-1}}h_L(\xi)d\mu.$$
Therefore,  the following optimization problem is naturally considered
$$m=\sup\limits_{K\in\mathcal{K}}\{F(K^\circ): W_\mu(K)=1\}$$
in our proof of Theorem 1.2. This is in some sense different from the usual setting. As a byproduct of our setting, we have the following isoperimetric inequality
\vskip 0.05in
{\em{\bf Theorem 1.3:} Let $H^{N-1}$ denote the $(N-1)$-dimensional Hausdorff measure and $B$ denote a ball in $\mathbb{R}^N$. Then, for any bounded convex domain $\Omega$, there holds $F(\Omega)\leq F(B)$ provided that $W_{H^{N-1}}(\Omega)=W_{H^{N-1}}(B)$. Moreover $F(\Omega)=F(B)$ and $W_{H^{N-1}}(\Omega)=W_{H^{N-1}}(B)$ imply $\Omega=B$ up to a translation.} 
\vskip 0.05in
{\em{\bf Remark 1.4:} Since $W_{H^{N-1}}(B_1(0))=2$ and $F$ is homogeneous of oder $\alpha=N+\frac{2(1+\beta)}{1-\beta}$, Theorem 1.3 can be rewritten as
$$F(\Omega)\leq\frac{F(B_1)}{2^\alpha}W^\alpha_{H^{N-1}}(\Omega)$$
for any bounded convex domain $\Omega$, and the equality if and only if $\Omega$, up to a translation, is a ball.}
\vskip 0.1in

The rest of this paper is arranged in the following way. Section 2 contains some preliminary results. Section 3 devotes to prove Theorem 1.1. Section 4 contains a proof of variational formula for $F(\Omega)$ in general class of bounded convex domains. The proof of Theorem 1.2 is contained in section 5.
\section*{2. Preliminary Results}

\setcounter{section}{2}

\setcounter{equation}{0}
Let $\Omega$ be a bounded domain in $\mathbb{R}^N$ with Lipschitz boundary $\partial\Omega$ and $H_0^1(\Omega)$ be the standard Soblev space with norm $\|u\|_{H_0^1}=(\int_\Omega|\nabla u|^2dx)^{\frac{1}{2}}$. For $\beta\in (0,1)$, we consider the sub-linear problem
\begin{equation}\label{eq221}
\left\{\begin{array}{ll}
-\Delta\varphi=\varphi^\beta & x\in \Omega\\
\varphi>0 &x\in\Omega\\
\varphi=0 & x\in\partial\Omega.
\end{array}
\right.
\end{equation}
It is well known that problem (\ref{eq221}) has a unique solution $\varphi(x)\in H_0^1(\Omega)$ which can be obtained as a minimizer of functional
$$I(u)=\frac{1}{2}\int_\Omega|\nabla u|^2dx-\frac{1}{\beta+1}\int_\Omega |u|^{\beta+1}dx$$
defined on $H_0^1(\Omega)$. Therefore, $\varphi(x)$ verifies
\begin{equation}\label{eq220}
\int_\Omega\nabla\varphi\cdot\nabla vdx-\int_\Omega\varphi^\beta vdx=0\ \ \ \ \forall\ v\in H_0^1(\Omega).
\end{equation}
Consequently
\begin{equation}\label{eq222}
\int_\Omega|\nabla\varphi|^2dx=\int_\Omega\varphi^{\beta+1}dx.
\end{equation}
Moreover, by the standard regularity theory of elliptic equations, we also have $\varphi(x)\in C^2(\Omega)\cap C(\overline{\Omega})$. If we assume further that $\Omega$ is convex, then the solution $\varphi(x)$ of problem (\ref{eq221}) possess the following important proposition.
\vskip 0.05in
{\em {\bf Proposition 2.1:} Let $\Omega$ be a bounded convex domain in $\mathbb{R}^N$ and $\varphi(x)$ be a solution to problem (\ref{eq221}). If we denote by $R_\Omega$ the diameter of $\Omega$, then there is a positive constant $C(N,R_\Omega,\beta)$ depending only on the dimension $N$, the diameter $R_\Omega$ of $\Omega$ and the number $\beta$ such that
$$\sup\limits_{x\in\Omega}|\varphi(x)|+\sup\limits_{x\in\Omega}|\nabla\varphi(x)|\leq C(N,R_\Omega,\beta).$$}
\vskip 0.05in
{\bf Proof:} By Theorem 3.7 of \cite{GT2001}, we have
$$\sup\limits_{x\in\Omega}|\varphi(x)|\leq(e^{2R_\Omega}-1)\sup\limits_{x\in\Omega}|\varphi(x)|^\beta.$$
From this and the fact that $0<\beta<1$, we obtain
$$\sup\limits_{x\in\Omega}|\varphi(x)|\leq(e^{2R_\Omega}-1)^{\frac{1}{1-\beta}}.$$

To estimate $\sup\limits_{x\in\Omega}|\nabla\varphi(x)|$, we recall a result of \cite{C2005} which says that $\varphi^{\frac{1-\beta}{2}}(x)$ is concave on $\Omega$. This implies that, for any $t\geq 0$, the level set $\Omega_t$ of $\varphi(x)$ at level $t$ defined by
$$\Omega_t=\{x\in\Omega:\varphi(x)>t\}$$
is convex. For any $x_0\in\Omega$, we let $t=\varphi(x_0)>0$ and $M_\Omega=\sup\limits_{x\in\Omega}\varphi(x)$. If $\varphi(x_0)=M_\Omega$, then $\nabla\varphi(x_0)=0$ and the claim is true. Assume $\varphi(x_0) < M_\Omega$, this implies that $x_0\in\partial\Omega_t$. The convex set $\Omega_t$ admits a support hyperplane $T_0$ at $x_0$. We may choose an orthogonal coordinate system with origin $O$ and coordinates $x_1, \cdots , x_n$ in $R^n$ such that $x_0=O$, $T_0=\{x\in R^n: x_n =0\}$ and $\Omega_t\subset\{x\in R^n: 0\leq x^n\leq R_\Omega\}$. By a standard argument based on the implicit function theorem, $\partial\Omega_t$ is of class $C^\infty$ so that $T_0$ is in fact the tangent hyperplane to $\partial\Omega_t$ at $x_0$. Consequently we have
$$|\nabla\varphi(x_0)| =\frac{\partial\varphi}{\partial x_n}(x_0).$$
Considering the auxiliary function
$$w(x) = w(x_1, \cdots, x_n) = t + \frac{1}{2}M^\beta_\Omega x_n(R_\Omega -x_n), \quad x\in R^n,$$
we have $-\Delta w(x) =M^\beta_\Omega\geq\varphi^\beta(x)=-\Delta\varphi(x)$ fore very $x\in\Omega_t$ and $w(x)\geq t=\varphi(x)$ for every $x\in\partial\Omega_t$. Therefore, it follows from the comparison Principle that
$$w(x)\geq\varphi(x)\quad\quad\quad\quad x\in\Omega_t.$$
Taking $\varphi(x_0) = w(x_0)$ into account, we finally have
$$\frac{\partial\varphi}{\partial x_n}( x_0) \leq\frac{\partial w}{\partial x_n}(x_0)=\frac{1}{2}M_\Omega^\beta R_\Omega.$$
Therefore
$$\sup\limits_{x\in\Omega}|\nabla\varphi(x)|\leq \frac{1}{2}(e^{2R_\Omega}-1)^{\frac{\beta}{1-\beta}}R_\Omega.$$
This completes the proof of Proposition 2.1.
\vskip 0.05in

To investigate the boundary behavior of $\nabla\varphi$, we briefly recall the notion of non-tangential limit of a function at a boundary point of a domain (for more details we refer the reader to \cite{K}). Let $\Omega$ be an open bounded subset of $\mathbb{R}^{N}$ and $d(y,\partial\Omega)$ be the distance between $y$ and $\partial\Omega$ which is defined by
$$d(y,\partial\Omega)=\inf\limits_{x\in\partial\Omega}d(y,x).$$

For $x\in\partial \Omega$ and $\alpha>0$, the non-tangential cone with vertex $x$ is defined by 
$$\Gamma_{\alpha}(x)=\{y\in \Omega:|y-x| \leq(1+\alpha)d(y,\partial\Omega)\}.$$
A sequence of points $y_{i} \in\Omega$ is said to converge non-tangentially to $x\in\partial \Omega$,  if $y_i\in\Gamma_{\alpha}(x)$ for some $\alpha>0$ and
$$\lim _{i \rightarrow\infty} y_{i}=x.$$
A function $w(y)$ defined in $\Omega$ is said to admit a non-tangential limit $L$ at $x\in\partial\Omega$ if for every $\alpha>0$ there holds
$$\lim _{y\rightarrow x, y\in\Gamma_{\alpha}(x)} w(y)=L$$
which is simply denoted by
$$\lim _{y\rightarrow x\text { n.t. }} w(x)=L.$$
For proceeding on, we need the following important result about harmonic functions by Dahlberg (see \cite{D}).
\vskip 0.05in
 {\em {\bf Theorem D:} Let $\Omega$ be a bounded Lipschitz domain in $\mathbb{R}^N$ and $w(x)$ be a harmonic in $\Omega$. If $w(x)$ is bounded below, then it has finite non-tangential limit for surface measure $\sigma$-almost every point $y\in\partial\Omega$.}
\vskip 0.05in
With the aid of Theorem D, we can prove the following proposition of $\nabla\varphi$.
\vskip 0.05in
{\em {\bf Proposition 2.2:} Let $\Omega$ be a bounded convex domain in $\mathbb{R}^N$ and $\varphi(x)$ be a solution to problem (\ref{eq221}). Then, each component of $\nabla\varphi(x)$ has finite non-tangential  limit for surface measure $\sigma$-almost every point $y\in\partial\Omega$.}
\vskip 0.05in
{\bf Proof:}\ Let $\Gamma(x-\xi)$ be the fundamental solution of Laplace equation. That is $\Gamma(x-\xi)$ is given by
$$
\Gamma(x-\xi)=\left\{\begin{array}{ll}
\frac{1}{N(2-N)\omega_N}|x-\xi|^{2-N}&N>2,\\
\frac{1}{2\pi}\log|x-\xi|&N=2.
\end{array}
\right.
$$
Where $\omega_N$ is the volume of unit ball in $\mathbb{R}^N$.

Let $v(x)$ be the Newton's potential of $\varphi^\beta(x)$. That is $v(x)$ is given by
$$v(x)=\int_\Omega\Gamma(x-\xi)\varphi^\beta(\xi)d\xi.$$
Since $\varphi^\beta(x)$ is bounded in $\Omega$, we have $v(x)\in C^1(\mathbb{R}^N)$. Moreover, for any $x\in\Omega$, $\nabla v(x)$ can be expressed by
$$\nabla v(x)=\int_\Omega\nabla_x\Gamma(x-\xi)\varphi^\beta(\xi)d\xi.$$
From this and Lemma 2.1, we can conclude that there exists a constant $C(N,R_\Omega,\beta)$ depending only on the dimension $N$, the diameter $R_\Omega$ of $\Omega$ and the number $\beta$ such that 
\begin{equation}\label{eq223}
\sup\limits_{x\in\Omega} |\nabla v(x)|\leq C(N,R_\Omega,\beta).
\end{equation}
Noting $\varphi(x)>0$ in $\Omega$ and $\varphi(x)\in C^2(\Omega)$, we see that $\varphi^\beta(x)$ is locally H\'older continuous in $\Omega$. Therefore, $v(x)\in C^2(\Omega)$ and satisfies
\begin{equation}\label{eq224}
\begin{array}{ll}
-\Delta v(x)=\varphi^\beta(x) & x\in\Omega.
\end{array}
\end{equation}
Let $w(x)=\varphi(x)-v(x)$. Then, it follows from (\ref{eq221}) and (\ref{eq224}) that $w(x)$ is a harmonic function in $\Omega$. By the property of harmonic function, we know that each component of $\nabla w(x)$ is also harmonic in $\Omega$. Furthermore, we can see from Lemma 2.1 and (\ref{eq223}) that each component of $\nabla w(x)$ is bounded from below. Since convex domain is Lipschitz, it follows from Theorem D that each component of $\nabla w(x)$ has a finite non-tangential limit for surface measure $\sigma$-almost every point $y\in\partial\Omega$. Noting $\nabla\varphi(x)=\nabla w(x)+\nabla v(x)$ for any $x\in\Omega$, we can conclude that each component of $\nabla\varphi(x)$ has a finite non-tangential limit for surface measure $\sigma$-almost every point $y\in\partial\Omega$.

By Proposition 2.2, we know that if $\Omega$ is a bounded convex domain in $\mathbb{R}^N$, then, for $\sigma$-almost every point $x\in\partial\Omega$, $\nabla\varphi(x)$ can be defined by
$$\nabla\varphi(x)=\lim\limits_{y\to x\text{ n.t.}}\nabla\varphi(y).$$ 
Therefore, we can introduce the following definition.
\vskip 0.05in
{\em{\bf Definition 2.3:} Let $\Omega$ be a bounded convex domain in $\mathbb{R}^N$ and $\varphi(x)$ be the unique solution of problem (\ref{eq221}). Setting $d\mu_\Omega=g_*(|\nabla\varphi(x)|^2d\sigma)$, or in other words
$$\int_{\mathbb{S}^{N-1}}f(\xi)d\mu_{\Omega}(\xi)=\int_{\partial\Omega}f(\nu(x))|\nabla\varphi(x)|^2d\sigma$$
for every $f(\xi)\in C(\mathbb{S}^{N-1})$. Then it is easy to see that $\mu_\Omega$ is a non-negative measure on $\mathbb{S}^{N-1}$.}
\vskip 0.05in
Next, we consider the following domain functional
$$F(\Omega)=\int_\Omega |\nabla\varphi(x)|^2dx$$
which is an important ingredient for solving Minkowski problem associated with $\mu_\Omega$ given in definition 2.3. To investigate the proposition of $F(\Omega)$, we prove the following result first.
\vskip 0.05in
{\em{\bf Proposition 2.4:} Let $\Omega_1$ and $\Omega_2$ be two bounded domains in $\mathbb{R}^N$ such that $\Omega_1\subset\Omega_2$. If $\varphi_1(x)$ and $\varphi_2(x)$ are respectively the solution of problem (\ref{eq221}) on $\Omega_1$ and $\Omega_2$, then $\varphi_1(x)\leq\varphi_2(x)$ for any $x\in\Omega_1$.}
\vskip 0.03in
{\bf Proof:} Let $A=\frac{1+\beta}{1-\beta}$, $\quad B=\frac{2}{1+\beta}$, and $q=\frac{2}{1-\beta}$. For a solution $\varphi(x)$ to problem (\ref{eq221}) in $\Omega$, we can easily see that $v(x)=\varphi^{\frac{1}{q}}(x)$ is a solution to the following problem
\begin{equation*}
\left\{\begin{array}{ll}
-\Delta v=\frac{1}{v}(A|\nabla v|^{2}+B)&x\in\Omega\\
v>0 & x\in\Omega \\
v=0 &x\in\partial\Omega
\end{array}
\right.
\end{equation*}
Setting $w(x)=v_{1}(x)-v_{2}(x)$, we can easily check that $w(x)$ satisfies
\begin{equation}\label{c11}
\left\{
\begin{array}{ll}
\Delta w(x)=\frac{1}{v_1(x)v_2(x)}(Bw(x)+A|\nabla v_1(x)|^2w(x)+Av_1(x)(|\nabla v_2(x)|^2-|\nabla v_1(x)|^2))& x\in\Omega_1\\
w(x)=-v_2(x)=-\varphi_2^{\frac{2}{1-\beta}(x)}<0&x\in\partial\Omega_1.
\end{array}
\right.
\end{equation}
To reach the conclusion of Proposition 2.4, we claim that
$$w(x)\leq 0\quad\mbox{for all}\ \ x\in\Omega_1.$$
Otherwise, there exits a point $x_{0}\in\Omega_{1}$ such that $$w(x_{0})=\max_{x\in\Omega_{1}}w(x)>0$$. Therefore, we have
$$\Delta w(x_{0})\leq 0,\quad \nabla w(x_0)=0,\quad |\nabla v_1(x_0)|=|\nabla v_2(x_0)|.$$
On the other hand, by the first equation in (\ref{c11}), we have
$$
\Delta w(x_0)=\frac{Bw(x_0)+A|\nabla v_1(x_0)|^2w(x_0)+Av_1(x_0)(|\nabla v_2(x_0)|^2-|\nabla v_1(x_0)|^2)}{v_1(x_0)v_2(x_0)}=\frac{Bw(x_0)+A|\nabla v_1(x_0)|^2w(x_0)}{v_1(x_0)v_2(x_0)}>0.
$$
This is a contradiction with $\Delta w(x_0)\leq 0$. Therefore, $w(x)\leq 0$ for all $x\in\Omega_{1}$. This completes the proof of Proposition 2.4.
\vskip 0.05in
With the aid of Proposition 2.4, we can prove the following result.
\vskip 0.05in
{\em{\bf Proposition 2.5:} Let $\Omega$ be a bounded domain in $\mathbb{R}^N$. For any fixed $x_0\in\mathbb{R}^N$ and any $t>0$, we denote by $x_0+\Omega$ the set $\{x_0+x: x\in\Omega\}$ and by $t\Omega$ the set $\{tx: x\in\Omega\}$. Then the following statements are true.
\vskip 0.03in
(i)\ $F(x_0+\Omega)=F(\Omega)$.
\vskip 0.03in
(ii)\ $\mu_{t\Omega}=t^{\alpha-1}\mu_\Omega$ and $F(t\Omega)=t^\alpha F(\Omega)$ with $\alpha=N+\frac{2(1+\beta)}{1-\beta}$.
\vskip 0.03in
(iii)\ If $\Omega_1\subset\Omega_2$, then $F(\Omega_1)\leq F(\Omega_2)$.}
\vskip 0.05in
{\bf Proof:} The conclusion $(i)$ follows from the fact that the form of problem (\ref{eq221}) is translation invariant. The conclusion $(ii)$ follows from the fact that if $\varphi(x)$ is a solution of problem (\ref{eq221}) on $\Omega$, then $v(y)=t^{\frac{2}{1-\beta}}u(y/t)$ is a solution of problem (\ref{eq221}) on $t\Omega$. The conclusion $(iii)$ follows from Proposition 2.4 and equality (\ref{eq222}). In fact, If we denote by $\varphi_1(x)$ and $\varphi_2(x)$ respectively the solution of problem (\ref{eq221}) on $\Omega_1$ and $\Omega_2$, then it follows from $\Omega_1\subset\Omega_2$ and Proposition 2.4 that $\varphi_1(x)\leq\varphi_2(x)$ for any $x\in\Omega_1$. Taking (\ref{eq222}) into account,
we have
$$F(\Omega_1)=\int_{\Omega_1}|\nabla\varphi_1(x)|^2dx=\int_{\Omega_1}\varphi_1^{\beta+1}(x)dx\leq\int_{\Omega_2}\varphi_2^{\beta+1}(x)dx=\int_{\Omega_2}|\nabla\varphi_2(x)|^2dx=F(\Omega_2).$$
This completes the proof of Proposition 2.5
\vskip 0.05in
The following result is crucial for our proving of Theorem 1.1 and Theorem 1.2.
\vskip 0.05in
{\em{\bf Proposition 2.6:} Let $\Omega$ be a bounded convex domain and $\{\Omega_n\}_{n=1}^\infty$ be a sequence of bounded domains such that $\Omega_n\subset\Omega$ for arbitrary $n$ and $\lim\limits_{n\to\infty}d(\overline{\Omega}_n,\overline{\Omega})=0$. Denote by $\varphi(x)$ and $\varphi_n(x)$ respectively the solution of problem (\ref{eq221}) on $\Omega$ and $\Omega_n$. If the extension $\widetilde{\varphi}_n(x)$ of $\varphi_n(x)$ on $\Omega$ is defined by
$$
\widetilde{\varphi}_n(x)=\left\{\begin{array}{ll}
\varphi_n(x)&x\in\Omega_n\\
0 &x\in\Omega\setminus\Omega_n.
\end{array}
\right.
$$
then the following statements are true.
\vskip 0.03in
(i) $\lim\limits_{n\to\infty}\widetilde{\varphi}_n(x)=\varphi(x)$ in $\Omega$.
\vskip 0.03in
(ii) $\lim\limits_{n\to\infty}F(\Omega_n)=F(\Omega)$.
\vskip 0.03in
(iii) $\lim\limits_{n\to\infty}\int_{\Omega_n}|\nabla\varphi_n(x)-\nabla\varphi(x)|^2dx=0$.}
\vskip 0.05in
{\bf Proof:} Since the form of problem (\ref{eq221}) is translation invariant and $\Omega$ is an open convex set, we may properly choose origin so that 
$B_\rho=\{x\in\mathbb{R}^N: |x|<\rho\}\subset\Omega$ for some $\rho>0$. Therefore, the support function $h_{\overline{\Omega}}(\xi)$ has a uniform lower bound $\rho$. This implies that
$$|\frac{h_{\overline{\Omega}_n} (\xi)}{h_{\overline{\Omega}} (\xi)}-1|\leq\frac{1}{\rho} |h_{\overline{\Omega}_n}(\xi)-h_{\overline{\Omega}}(\xi)|\leq\frac{1}{\rho} d(\overline{\Omega}_n,\overline{\Omega})$$
and $\frac{h_{\overline{\Omega}_n(\xi)}}{h_{\overline{\Omega}} (\xi)}$ converges uniformly to $1$ on $\mathbb{S}^{N-1}$. Consequently, for any $\varepsilon>0$, there is a positive integer $N(\varepsilon)$ depending only on $\varepsilon$ such that
$$h_{\overline{\Omega}_n}(\xi)>(1-\varepsilon)h_{\overline{\Omega}} (\xi)\quad\forall n>N(\varepsilon),\ \ \xi\in\mathbb{S}^{N-1}.$$
Therefore,
$$(1-\varepsilon)\Omega\subset\Omega_n\subset\Omega\quad\forall n>N(\varepsilon).$$
It is easy to check that $v(y)=(1-\varepsilon)^{\frac{2}{1-\beta}}\varphi(\frac{y}{1-\varepsilon})$ is a solution of problem (\ref{eq221}) in $(1-\varepsilon)\Omega$. By Proposition 2.4, we have $v(y)\leq\varphi_n(y)\leq\varphi(y)$ for any $y\in (1-\varepsilon)\Omega$. Denoting by $\widetilde{v}(y)$ the extension of $v(y)$ on $\Omega$ defined by
$$\widetilde{v}(y)=\left\{\begin{array}{ll}
v(y)&y\in (1-\varepsilon)\Omega\\
0 &\Omega\setminus (1-\varepsilon)\Omega,
\end{array}
\right.
$$
we can see that
$$\widetilde{v}(y)\leq\widetilde{\varphi}_n(y)\leq\varphi(y)\quad\forall n>N(\varepsilon),\ \ y\in\Omega.$$
From this, we can conclude that, for any fixed $y\in\Omega$, there holds
$$\widetilde{v}(y)\leq\varliminf\limits_{n\to\infty}\widetilde{\varphi}_n(y)\leq\varlimsup\limits_{n\to\infty}\widetilde{\varphi}_n(y)\leq\varphi(y).$$
Sending $\varepsilon$ to zero in the above inequality and noting $\lim\limits_{\varepsilon\to 0}\widetilde{v}(y)=\varphi(y)$, we get
$$\lim\limits_{n\to\infty}\widetilde{\varphi}_n(y)=\varphi(y).$$
This proves the statement $(i)$.

 By $(ii)$ and $(iii)$ of proposition 2.5, we have
$$(1-\varepsilon)^\alpha F(\Omega)\leq F(\Omega_n)\leq F(\Omega)\quad\forall n>N(\varepsilon).$$
This leads to the following inequality
$$(1-\varepsilon)^\alpha F(\Omega)\leq \varliminf\limits_{n\to\infty}F(\Omega_n)\leq\varlimsup\limits_{n\to\infty}F(\Omega_n)\leq F(\Omega).$$
Sending $\varepsilon$ to zero in the above inequality, we obtain
$$\lim\limits_{n\to\infty}F(\Omega_n)=F(\Omega).$$
This proves the statement $(ii)$.

By Proposition 2.1 and $\varphi_n(x)\in H_0^1(\Omega_n)$, we can easily see that $\widetilde{\varphi}_n(x)\in H_0^1(\Omega)$ and
$$\nabla\widetilde{\varphi}_n(x)=\left\{\begin{array}{ll}
\nabla\varphi_n(x)&x\in\overline{\Omega}_n\\
0&x\in\Omega\setminus\overline{\Omega}_n.
\end{array}
\right.
$$
Taking $\widetilde{\varphi}_n(y)$ as a test function in (\ref{eq220}), we get
$$\int_\Omega\widetilde{\varphi}_n(x)\varphi^\beta(x)dx=\int_\Omega\nabla\widetilde{\varphi}_n\cdot\nabla\varphi dx=\int_{\Omega_n}\nabla\varphi_n\cdot\nabla\varphi dx.$$
By statement $(i)$, Proposition 2.1 and the dominated convergence theorem, we reach
$$\lim\limits_{n\to\infty}\int_\Omega\widetilde{\varphi}_n(x)\varphi^\beta(x)dx=\int_\Omega\varphi^{\beta+1}(x)dx=\int_\Omega|\nabla\varphi(x)|^2dx.$$
Therefore
\begin{equation}\label{eq225}
\lim\limits_{n\to\infty}\int_{\Omega_n}\nabla\varphi_n\cdot\nabla\varphi dx=\int_\Omega|\nabla\varphi(x)|^2dx.
\end{equation}
Since
$$\int_{\Omega_n}|\nabla\varphi_n(x)-\nabla\varphi(x)|^2dx=F(\Omega_n)-2\int_{\Omega_n}\nabla\varphi_n\cdot\nabla\varphi dx+\int_{\Omega_n}|\nabla\varphi(x)|^2dx,$$
it follows from statement $(ii)$, (\ref{eq225}) and the fact
$$\lim\limits_{n\to\infty}\int_{\Omega_n}|\nabla\varphi(x)|^2dx=F(\Omega)$$
that
$$\lim\limits_{n\to\infty}\int_{\Omega_n}|\nabla\varphi_n(x)-\nabla\varphi(x)|^2dx=0.$$
This proves the statement $(iii)$ and thus completes the proof of Proposition 2.6.

\section*{3. Weak continuity of $\mu_\Omega$}

\setcounter{section}{3}

\setcounter{equation}{0}

\noindent
This section devotes to prove the weak continuity of measure $\mu_\Omega$ defined in definition 2.3 with respect to the Hausdorff metric. Our main aim is to prove the following theorem.
\vskip 0.05in
{\em{\bf Theorem 3.1:} Let $\Omega$ and $\Omega_n$ be bounded convex domains in $\mathbb{R}^N$. If $d_H(\overline{\Omega}_n,\overline{\Omega})\to 0$ as $n\to\infty$, then $\mu_{\Omega_n}$ converges weakly to $\mu_\Omega$ in the sense of measure, that is $\lim\limits_{n\to\infty}\int_{\mathbb{S}^{N-1}}f(\xi)d\mu_{\Omega_n}=\int_{\mathbb{S}^{N-1}}f(\xi)d\mu_{\Omega}$ for any $f\in C(\mathbb{S}^{N-1})$.}
\vskip 0.05in
To reach our aim, we quote here a results of A. Colesanti and M. Fimiani  given in \cite{CF} as our Lemma 3.2. 
\vskip 0.05in
{\em{\bf Lemma 3.2:} Let $\Omega$ and $\Omega_n$ be bounded convex domains of $\mathbb{R}^N$ such that $\Omega_n\subset\Omega$ for any $n$ and $d_H(\overline{\Omega}_n,\overline{\Omega})\to 0$ as $n\to\infty$. Let $f: \partial\Omega\to\mathbb{R}$ and $f_n: \partial\Omega_n\to\mathbb{R}$ be $\sigma$-measurable functions. If $f$ and $f_n$ satisfy
\vskip 0.03in
(i) There exists a universal constant $C>0$ such that $\|f\|_{L^\infty(\partial\Omega)}\leq C$ and $\|f_n\|_{L^\infty(\partial\Omega)}\leq C$ for any $n$,
\vskip 0.03in
(ii) For $\sigma$-almost every $x\in\partial\Omega$, if $x_n\in\partial\Omega_n$, $f_n$ is well-defined at $x_n$ and $\lim\limits_{n\to\infty\ \text{n.t.}}x_n=x$ imply 
           $$\lim\limits_{n\to\infty}f_n(x_n)=f(x),$$
then there holds
$$\lim\limits_{n\to\infty}\int_{\partial\Omega_n}f_n(x)d\sigma=\int_{\partial\Omega}f(x)d\sigma.$$}
\vskip 0.05in
{\em{\bf Lemma 3.3:} Let $\Omega$ and $\Omega_n$ be bounded convex domains of $\mathbb{R}^N$ such that $\Omega_n\subset\Omega$ for any $n$ and $d_H(\overline{\Omega}_n,\overline{\Omega})\to 0$ as $n\to\infty$. Let $\varphi(x)$ be the unique solution to problem (\ref{eq221}) in $\Omega$. Let $x\in\partial\Omega$ be such that at which $\nabla\varphi$ has finite non-tangential limit and $\partial\Omega$ is differentiable. Let $x_n\in\Omega_n$ be such that $\partial\Omega_n$ is differentiable at $x_n$ and $\lim\limits_{n\to\infty\ \text{n.t.}}x_n=x$. Define $\nabla_{T_n}\varphi(x_n)$ by
$$\nabla_{T_n}\varphi(x_n)=\nabla\varphi(x_n)-\nu_n(x_n)\cdot\nabla\varphi(x_n)\nu_n(x_n)$$
with $\nu_n(x_n)$ being the outer unit vector normal to $\partial\Omega_n$ at $x_n$. Then $\lim\limits_{n\to\infty}\nabla_{T_n}\varphi(x_n)=0$.}
\vskip 0.05in
{\bf Proof:} If $\lim\limits_{y\to x\ \text{n.t.}}\nabla\varphi(y)=0$, then we are done. Therefore, we assume that $\lim\limits_{y\to x\ \text{n.t.}}\nabla\varphi(y)\neq 0$. This implies that 
$|\nabla\varphi(x_n)|>0$ for large enough $n$. Let $\Omega_{\varepsilon_n}=\{y\in\Omega: u(y)>\varepsilon_n\}$ with $\varepsilon_n=u(x_n)>0$ and $\widetilde{\nu}_n(y)$ be the outer unit vector normal to $\partial\Omega_{\varepsilon_n}$ at $y$. Since $|\nabla\varphi(x_n)|>0$, $\frac{\nabla\varphi(x_n)}{|\nabla\varphi(x_n)|}$ is well-defined for large enough $n$ and 
$$-\frac{\nabla\varphi(x_n)}{|\nabla\varphi(x_n)|}=\widetilde{\nu}_n(x_n)$$
due to $u(y)=\varepsilon_n$ for any $y\in\Omega_{\varepsilon_n}$. Noting $\lim\limits_{n\to\infty}\varepsilon_n=0$, we can see that $d_H(\overline{\Omega}_{\varepsilon_n},\overline{\Omega})\to 0$ as $n\to\infty$. Therefore $\lim\limits_{n\to\infty}\widetilde{\nu}_n(x_n)=\nu(x)$. By the assumption that $d_H(\overline{\Omega}_n,\overline{\Omega})\to 0$ as $n\to\infty$, we also have $\lim\limits_{n\to\infty}\nu_n(x_n)=\nu(x)$. Rewritten $\nabla_{T_n}\varphi(x_n)$ in the following way
$$\nabla_{T_n}\varphi(x_n)=-|\nabla\varphi(x_n)|\widetilde{\nu}_n(x_n)+|\nabla\varphi(x_n)|(\widetilde{\nu}_n(x_n)\cdot\nu_n(x_n))\nu_n(x_n),$$
we can easily see that
$$\lim\limits_{n\to\infty}\nabla_{T_n}\varphi(x_n)=0.$$
\vskip 0.05in
The following lemma is a main step in our proof of Theorem 3.1.
\vskip 0.05in
{\em{\bf Lemma 3.4:} Let $\Omega$ and $\Omega_n$ be convex domains such that $\Omega_n\subset\Omega$ for any $n$ and $\Omega_n$ converges to $\Omega$ with respect to Hausdorff metric. Denote by $\varphi(x)$ and $\varphi_n(x)$ respectively the solution of problem (\ref{eq221}) in $\Omega$ and $\Omega_n$. Then we have
$$\lim\limits_{n\to\infty}\int_{\partial\Omega_n}|\nabla\varphi_n(x)-\nabla\varphi(x)|^2d\sigma=0.$$}
\vskip 0.04in
{\bf Proof:} Because the form of problem (\ref{eq221}) is translation invariant, we may choose origin $o$ properly so that the ball $B_{2\rho}(o)\subset\Omega$ for some fixed number $\rho>0$. Therefore, the support function $h_{\overline{\Omega}}$ of $\overline{\Omega}$ verifies $h_{\overline{\Omega}}(\xi)>2\rho$ for any $\xi\in\mathbb{S}^{N-1}$. Noting the support function $h_{\overline{\Omega}_n}$ of $\overline{\Omega}_n$ converges uniformly to $h_{\overline{\Omega}}$, we can conclude that there is a positive integer $N(\rho)$ such that $h_{\overline{\Omega}_n}(\xi)>\frac{3}{2}\rho$ for any $\xi\in\mathbb{S}^{N-1}$. The main difficult in estimating of $\int_{\partial\Omega_n}|\nabla\varphi_n(x)-\nabla\varphi(x)|^2d\sigma$ is that the boundary behavior of $|\nabla\varphi_n(x)-\nabla\varphi(x)|$ is not so clear. Our strategy is try to transfer the estimation of $\int_{\partial\Omega_n}|\nabla\varphi_n(x)-\nabla\varphi(x)|^2d\sigma$ to the estimation of volume integration. The usual tool for this transfer is divergence theorem. However, $\partial\Omega_n$ is not smooth enough for using divergence theorem. So, we should make a minor modification. For any $\varepsilon>0$, we consider the domain
$$\Omega_n^{\varepsilon}=\{x\in\Omega_n: \varphi_(x)>\varepsilon\}.$$
Since $\varphi_n(x)$ is smooth in $\Omega_n$, $\partial\Omega_n$ is smooth enough for using divergence theorem. Noting $d_H(\overline{\Omega^\varepsilon}_n,\overline{\Omega}_n)\to 0$ as $\varepsilon\to 0$, and, for $\sigma$-almost every $x\in\partial\Omega_n$
$$\lim\limits_{y\to x\ \text{n.t.}}|\nabla\varphi_n(y)-\nabla\varphi(y)|^2=|\nabla\varphi_n(x)-\nabla\varphi(x)|^2,$$
it follows from Lemma 3.2 that 
$$\lim\limits_{\varepsilon\to 0}\int_{\partial\Omega_n^\varepsilon}|\nabla\varphi_n(x)-\nabla\varphi(x)|^2d\sigma
=\int_{\partial\Omega_n}|\nabla\varphi_n(x)-\nabla\varphi(x)|^2d\sigma.$$
For simplicity, we set $W_n(x)=\varphi_n(x)-\varphi(x)$ and $G_n(x)=\varphi_n^\beta(x)-\varphi^\beta(x)$. It follows from the equation satisfied by $\varphi_n(x)$ and $\varphi(x)$ that $W_n(x)$ verifies
\begin{equation}\label{eq331}
\left\{\begin{array}{ll}
-\Delta W_n=G_n&x\in\Omega_n^\varepsilon\\
\varphi_n=\varepsilon &x\in\partial\Omega_n^\varepsilon.
\end{array}
\right.
\end{equation}
Our next aim is to estimate $\int_{\partial\Omega_n^\varepsilon}|\nabla W_n|^2dx$. To reach our goal, we do the following computation first. By an easy computation, we have
\begin{equation}\label{eq332}
div((x\cdot\nabla W_n)\nabla W_n)=(x\cdot\nabla W_n)\Delta W_n+\nabla W_n\cdot\nabla(x\cdot\nabla W_n).
\end{equation}
\begin{equation}\label{eq333}
\nabla W_n\cdot\nabla(x\cdot\nabla W_n)=|\nabla W_n|^2+\sum\limits_{i,j=1}^Nx_j\frac{\partial W_n}{\partial x_i}\frac{\partial^2W_n}{\partial x_i\partial x_j}.
\end{equation}
Let $\nu_n^\varepsilon(x)$ denote the unit outer normal of $\partial\Omega_n^\varepsilon$. By making use of integration by part, we can get
\begin{equation}\label{eq334}
\int_{\Omega_n^\varepsilon}\sum\limits_{i,j=1}^Nx_j\frac{\partial W_n}{\partial x_i}\frac{\partial^2W_n}{\partial x_i\partial x_j}dx=\frac{1}{2}\int_{\partial\Omega_n^\varepsilon}(x\cdot\nu_n^\varepsilon(x))|\nabla W_n|^2d\sigma-\frac{N}{2}\int_{\Omega_n^\varepsilon}|\nabla W_n|^2dx.
\end{equation}
From (\ref{eq333}) and (\ref{eq334}), we have
\begin{equation}\label{eq335}
\int_{\Omega_n^\varepsilon}\nabla W_n\cdot\nabla(x\cdot\nabla W_n)dx=\frac{1}{2}\int_{\partial\Omega_n^\varepsilon}(x\cdot\nu_n^\varepsilon(x))|\nabla W_n|^2d\sigma-\frac{N-2}{2}\int_{\Omega_n^\varepsilon}|\nabla W_n|^2dx.
\end{equation}
Integrating the identity (\ref{eq332}) on $\Omega_n^\varepsilon$, we obtain
\begin{equation}\label{eq336}
\int_{\partial\Omega_n^\varepsilon}(x\cdot\nabla W_n)(\nu_n^\varepsilon\cdot\nabla W_n)d\sigma=-\int_{\Omega_n^\varepsilon}(x\cdot\nabla W_n)G_n(x) dx+\int_{\Omega_n^\varepsilon}\nabla W_n\cdot\nabla(x\cdot\nabla W_n)dx.
\end{equation}
Let $x_{T_n^\varepsilon}=x-(x\cdot\nu_n^\varepsilon(x))\nu_n^\varepsilon$. Geometrically, $x_{T_n^\varepsilon}$ is the projection of $x$ in the plane tangent to $\partial\Omega_n^\varepsilon$. Let $\nabla_{T_n^\varepsilon}W_n$ denote the tangential part of $\nabla W_n$ defined by $\nabla_{T_n^\varepsilon}W_n=\nabla W_n-(\nu_n^\varepsilon\cdot\nabla W_n)\nu_n^\varepsilon$. Then, we can easily check
\begin{equation}\label{eq337}
(x\cdot\nabla W_n)(\nu_n^\varepsilon\cdot\nabla W_n)=(x_{T_n^\varepsilon}\cdot\nabla_{T_n^\varepsilon}W_n)(\nu_n^\varepsilon\cdot\nabla W_n)+(x\cdot\nu_n^\varepsilon(x))(\nu_n^\varepsilon\cdot\nabla W_n)^2
\end{equation}
\begin{equation}\label{eq338}
(x\cdot\nu_n^\varepsilon)|\nabla W_n|^2=(x\cdot\nu_n^\varepsilon)|\nabla_{T_n^\varepsilon}W_n|^2+(x\cdot\nu_n^\varepsilon)(\nu_n^\varepsilon\cdot\nabla W_n)^2.
\end{equation}
From (\ref{eq335}), (\ref{eq336}), (\ref{eq337}) and (\ref{eq338}), we obtain
\begin{equation}\label{eq339}
\left\{\begin{array}{ll}
\int_{\partial\Omega_n^\varepsilon}(x\cdot\nu_n^\varepsilon)(\nu_n^\varepsilon\cdot\nabla W_n)^2d\sigma=I_1+I_2\\
I_1=\int_{\partial\Omega_n^\varepsilon}(x\cdot\nu_n^\varepsilon)|\nabla_{T_n^\varepsilon}W_n|^2d\sigma-2\int_{\partial\Omega_n^\varepsilon}(x_{T_n^\varepsilon}\cdot\nabla_{T_n^\varepsilon}W_n)(\nu_n^\varepsilon\cdot\nabla W_n)d\sigma\\
I_2=-2\int_{\Omega_n^\varepsilon}(x\cdot\nabla W_n)G_n(x)dx+(2-N)\int_{\Omega_n^\varepsilon}|\nabla W_n|^2dx.
\end{array}
\right.
\end{equation}
Since $\nu_n^\varepsilon$ converges uniformly to $\nu_n$ as $\varepsilon\to 0$ and $h_{\Omega_n}>\frac{3}{2}\rho$ for $n>N(\rho)$, we can see that $(x\cdot\nu_n^\varepsilon(x))>\rho$ on $\partial\Omega_n^\varepsilon$ for all $\varepsilon$ small enough and $n>N(\rho)$. Therefore
\begin{equation}\label{eq3310}
\int_{\partial\Omega_n^\varepsilon}(\nu_n^\varepsilon\cdot\nabla W_n)^2d\sigma\leq\frac{|I_1|+|I_2|}{\rho}.
\end{equation}
From this we have
\begin{equation}\label{eq3311}
\begin{array}{ll}
\int_{\partial\Omega_n^\varepsilon}|\nabla W_n|^2d\sigma &=\int_{\partial\Omega_n^\varepsilon}|\nabla_{T_n^\varepsilon} W_n|^2d\sigma+\int_{\partial\Omega_n^\varepsilon}(\nu_n^\varepsilon\cdot\nabla W_n)^2d\sigma\\
&\leq\int_{\partial\Omega_n^\varepsilon}|\nabla_{T_n^\varepsilon} W_n|^2d\sigma+\frac{|I_1|+|I_2|}{\rho}.
\end{array}
\end{equation}
Noting $\varphi_n(x)=\varepsilon$ for any $x\in\partial\Omega_n^\varepsilon$, we have $\nabla_{T_n^\varepsilon}\varphi_n=0$ and $\nabla_{T_n^\varepsilon}W_n=-\nabla_{T_n^\varepsilon}\varphi$. Therefore
\begin{equation}\label{eq3312}
I_1=\int_{\partial\Omega_n^\varepsilon}(x\cdot\nu_n^\varepsilon)|\nabla_{T_n^\varepsilon}\varphi|^2d\sigma-2\int_{\partial\Omega_n^\varepsilon}(x_{T_n^\varepsilon}\cdot\nabla_{T_n^\varepsilon}\varphi)(\nu_n^\varepsilon\cdot\nabla W_n)d\sigma
\end{equation}
\begin{equation}\label{eq3313}
\int_{\partial\Omega_n^\varepsilon}|\nabla W_n|^2d\sigma\leq\int_{\partial\Omega_n^\varepsilon}|\nabla_{T_n^\varepsilon}\varphi|^2d\sigma+\frac{|I_1|+|I_2|}{\rho}.
\end{equation}
By Proposition 2.1, we can see that there exists a positive constant $C$ depending only on the diameter $R_\Omega$ of $\Omega$ such that
\begin{equation}\label{eq3314}
|I_1|\leq C(\int_{\partial\Omega_n^\varepsilon}|\nabla_{T_n^\varepsilon}\varphi|^2d\sigma+\int_{\partial\Omega_n^\varepsilon}|\nabla_{T_n^\varepsilon}\varphi(x)|d\sigma.)
\end{equation}
\begin{equation}\label{eq3315}
|I_2|\leq C(\int_{\Omega_n^\varepsilon}|G_n(x)|dx+\int_{\Omega_n^\varepsilon}|\nabla W_n|^2dx.)
\end{equation}
From (\ref{eq3313}), (\ref{eq3314}) and (\ref{eq3315}), we can conclude that there is a positive constant $C$ independent of $\varepsilon$ and $n$ such that
\begin{equation}\label{eq3316}
\int_{\partial\Omega_n^\varepsilon}|\nabla W_n|^2d\sigma\leq C(\int_{\partial\Omega_n^\varepsilon}|\nabla_{T_n^\varepsilon}\varphi|^2d\sigma+\int_{\partial\Omega_n^\varepsilon}|\nabla_{T_n^\varepsilon}\varphi(x)|d\sigma+\int_{\Omega_n^\varepsilon}|G_n(x)|dx+\int_{\Omega_n^\varepsilon}|\nabla W_n|^2dx.)
\end{equation}
Because $\nabla_{T_n^\varepsilon}\varphi(x)$ converges to $\nabla_{T_n}\varphi(x)$ as $\varepsilon\to 0$. Sending $\varepsilon$ to zero in (\ref{eq3316}), we can infer from Lemma 3.2 that
\begin{equation}\label{eq3317}
\int_{\partial\Omega_n}|\nabla W_n|^2d\sigma\leq C(\int_{\partial\Omega_n}|\nabla_{T_n}\varphi|^2d\sigma+\int_{\partial\Omega_n}|\nabla_{T_n}\varphi(x)|d\sigma+\int_{\Omega_n}|G_n(x)|dx+\int_{\Omega_n}|\nabla W_n|^2dx.)
\end{equation}
By Proposition 2.6 and Lemma 3.3, we can see that
$$\lim\limits_{n\to\infty}(\int_{\Omega_n}|G_n(x)|dx+\int_{\Omega_n}|\nabla W_n|^2dx)=0.$$
$$\lim\limits_{n\to\infty}(\int_{\partial\Omega_n}|\nabla_{T_n}\varphi|^2d\sigma+\int_{\partial\Omega_n}|\nabla_{T_n}\varphi(x)|d\sigma)=0.$$
Therefore 
$$\lim\limits_{n\to\infty}\int_{\partial\Omega_n}|\nabla\varphi_n(x)-\nabla\varphi(x)|^2d\sigma=\lim\limits_{n\to\infty}\int_{\partial\Omega_n}|\nabla W_n|^2d\sigma=0.$$
This completes the proof of Lemma 3.4.
\vskip 0.05in
{\bf Proof of Theorem 3.1:} Let $\Omega$ and $\Omega_n$ be bounded domains of $\mathbb{R}^N$ such that $d_H(\overline{\Omega}_n,\overline{\Omega})\to 0$ as $n\to\infty$. Let $\nu$ and $\nu_n$ be the Gauss mapping of $\Omega$ and $\Omega_n$ respectively. Then $\nu_n$ converges uniformly to $\nu$. Let $\varphi(x)$ and $\varphi_n(x)$ be respectively the solution of problem (\ref{eq221}) in $\Omega$ and $\Omega_n$.  Denote by $\mu_\Omega$ the associated with $\Omega$. What we should do is proving
$$\lim\limits_{n\to\infty}\int_{\mathbb{S}^{N-1}}g(\xi)d\mu_{\Omega_n}(\xi)=\lim\limits_{n\to\infty}\int_{\mathbb{S}^{N-1}}g(\xi)d\mu_{\Omega}(\xi).$$
for any $g\in C(\mathbb{S}^{N-1})$. No loss of generality, we may assume that the origin $o$ belongs to $\Omega$ and that there is a positive number $\rho>0$ such that the ball $B_\rho(o)$ contains in $\Omega$. Therefore, the support function $h_{\overline{\Omega}}$ of $\Omega$ has a lower bound $\rho$. By the uniform convergence of the support function $h_{\overline{\Omega}_n}$ of $\Omega_n$ to $h_{\overline{\Omega}}$, we may also assume that $h_{\overline{\Omega}_n}$ has a lower bound $\frac{1}{2}\rho$ for all $n$. To make our previous results in use, we assume $\Omega_n\subset\Omega$ for any $n$. By the definition of $\mu_\Omega$, we have
\begin{equation}\label{eq3318}
\int_{\mathbb{S}^{N-1}}g(\xi)d\mu_{\Omega_n}(\xi)-\int_{\mathbb{S}^{N-1}}g(\xi)d\mu_{\Omega}(\xi)=\int_{\partial\Omega_n}g(\nu_n(x))|\nabla\varphi_n(x)|^2d\sigma-\int_{\partial\Omega}g(\nu(x))|\nabla\varphi(x)|^2d\sigma.
\end{equation}
Let $I_1$ and $I_2$ be defined by the following formula
\begin{equation}\label{eq3319}
I_1=\int_{\partial\Omega_n}g(\nu_n(x))(|\nabla\varphi_n(x)|^2-|\nabla\varphi(x)|^2)d\sigma.
\end{equation}
\begin{equation}\label{eq3320}
I_2=\int_{\partial\Omega_n}g(\nu_n(x))|\nabla\varphi(x)|^2d\sigma-\int_{\partial\Omega}g(\nu(x))|\nabla\varphi(x)|^2d\sigma.
\end{equation}
Then, (\ref{eq3318}) can be rewritten as
\begin{equation}\label{eq3321}
\int_{\mathbb{S}^{N-1}}g(\xi)d\mu_{\Omega_n}(\xi)-\int_{\mathbb{S}^{N-1}}g(\xi)d\mu_{\Omega}(\xi)=I_1+I_2.
\end{equation}
By the continuity of $g$ on compact set $\mathbb{S}^{N-1}$, we see that there exists a positive constant $C$ such that
$$|I_1|\leq C\int_{\partial\Omega_n}||\nabla\varphi_n(x)|^2-|\nabla\varphi(x)|^2|d\sigma.$$
By Proposition 2.1, we can conclude that 
$$\int_{\partial\Omega_n}||\nabla\varphi_n(x)|^2-|\nabla\varphi(x)|^2|d\sigma\leq C(\int_{\partial\Omega_n}|\nabla\varphi_n(x)-\nabla\varphi(x)|^2|d\sigma)^\frac{1}{2}$$
for some positive constant $C$.
Therefore, it follows from Lemma 3.4 that
$$\lim\limits_{n\to\infty}|I_1|\leq C\lim\limits_{n\to\infty}(\int_{\partial\Omega_n}|\nabla\varphi_n(x)-\nabla\varphi(x)|^2|d\sigma)^\frac{1}{2}=0.$$
Let $f_n(x)$ and $f((x)$ be defined by
$$f_n(x)=g(\nu_n(x))|\nabla\varphi(x)|^2,\quad f(x)=g(\nu(x))|\nabla\varphi(x)|^2.$$
Then, Proposition 2.1 and the continuity of $g$, can be used to ensure that there exists a positive constant $C$ such that
$$\|f_n\|_{L^\infty(\partial\Omega_n)}\leq C,\quad \|f\|_{L^\infty(\partial\Omega})\leq C.$$
Moreover, by the definition of $\nabla\varphi$ on $\partial\Omega$, the continuity of $g$ and the uniform convergence of $\nu_n$ to $\nu$, we further know that, for $\sigma$-almost every $x\in\partial\Omega$ and $x_n\in\partial\Omega_n$ there holds
$$\lim\limits_{n\to\infty}f_n(x_n)=f(x)$$
whenever $\lim\limits_{n\to\infty}x_n=x$. Therefore, Lemma 3.2 can used to imply that
$$\lim\limits_{n\to\infty}I_2=\lim\limits_{n\to\infty}(\int_{\partial\Omega_n}g(\nu_n(x))|\nabla\varphi(x)|^2d\sigma-\int_{\partial\Omega}g(\nu(x))|\nabla\varphi(x)|^2d\sigma)=0.$$
Consequently
$$\lim\limits_{n\to\infty}\int_{\mathbb{S}^{N-1}}g(\xi)d\mu_{\Omega_n}(\xi)=\int_{\mathbb{S}^{N-1}}g(\xi)d\mu_{\Omega}(\xi).$$
This completes the proof of Theorem 3.1 in the case $\Omega_n\subset\Omega$.

Now we assume that the condition $\Omega_n\subset\Omega$ is not fulfilled. By the uniform convergence of $h_{\overline{\Omega}_n}$ to $h_{\overline{\Omega}}$, we can fined a sequence of number $\alpha_n>0$ such that
$$\lim\limits_{n\to\infty}\alpha_n=1\quad\mbox{and}\quad \alpha_n\Omega_n\subset\Omega.$$
In fact, if choose $\alpha_n$ as
$$\alpha_n=\inf\limits_{\xi\in\mathbb{S}^{N-1}}\frac{h_{\overline{\Omega}}(\xi)}{h_{\overline{\Omega}_n}(\xi)}-\frac{1}{n},$$
then $\alpha_n$ is well-defined due to $h_{\overline{\Omega}}$ has a lower bound $\frac{1}{2}\rho$ and $\lim\limits_{n\to\infty}\alpha_n=1$ due to the uniform convergence of $h_{\overline{\Omega}_n}$ to $h_{\overline{\Omega}}$. Moreover, $\alpha_nh_{\overline{\Omega}_n}<h_{\overline{\Omega}}$ and $\alpha_n\Omega_n\subset\Omega$. Obviously, $d_H(\overline{\alpha_n\Omega_n},\overline{\Omega})\to 0$ as $n\to\infty$. Therefore, by the first step, we conclude that 
$$\lim\limits_{n\to\infty}\mu_{\alpha_n\Omega_n}=\mu_\Omega$$
weakly in the sense of measure. Noting that $\mu_\Omega$ is $\gamma=N+\frac{3\beta+1}{1-\beta}$-homogeneous, we have
$$\mu_{\alpha_n\Omega_n}=\alpha_n^\gamma\mu_{\Omega_n}.$$
Taking $\lim\limits_{n\to\infty}\alpha_n=1$ into account, we have
$$\lim\limits_{n\to\infty}\mu_{\Omega_n}=\lim\limits_{n\to\infty}\alpha_n^\gamma\mu_{\Omega_n}=\lim\limits_{n\to\infty}\mu_{\alpha_n\Omega_n}=\mu_\Omega.$$
weakly in the sense of measure. This completes the proof of Theorem 3.1.

\section*{4. Variational formula of $F$ for convex domains}

\setcounter{section}{4}

\setcounter{equation}{0}

\noindent
Let $K, L\in\mathcal{K}^N$. The Minkowski addition $K+L$ of $K$ and $L$ is defined by
$$K+L=\{x+y: x\in K, y\in L\}.$$
Minkowski addition has an equivalent description of support function which says that for arbitrary $K, L, Q\in\mathcal{K}^N$ there holds
$$Q=K+L\iff h_Q=h_K+h_L.$$
For any $K\in\mathcal{K}^N$, we denote by $K^\circ$ the interior of $K$ and by $\varphi_K(x)$ the unique solution of problem (\ref{eq221}) in $K^\circ$. Let $\mu_{K^\circ}$ and $F(K^\circ)$ be the measure and the functional introduced in section 2. For convenience, we recall that $\mu_{K^\circ}$ and $F(K^\circ)$ are defined by
$$d\mu_{K^\circ}=g_*(|\nabla\varphi_K|^2d\sigma),\quad F(K^\circ)=\int_{K^\circ}|\nabla\varphi_K|^2dx.$$
Where $g: \partial K\to\mathbb{S}^{N-1}$ is the Gauss mapping of $K$ which is given by $g(x)=\nu_K(x)$ and $\nu_K(x)$ is the outer unit vector normal to $\partial K$ at $x$.

A variational formula of $F$ for the interior of smooth convex bodies was given in \cite{HSX2018} by Y. Huang, C. Song and L.Xu which can be restated as
\vskip 0.05in
{\em{\bf Theorem HSX:} Let $K, L\in\mathcal{K}^N$ be of class $C^2$. Setting $K_{tL}=K+tL$ for $t>0$, then there holds
$$\frac{d}{dt}F(K_{tL}^\circ)|_{t=0}=\frac{1+\beta}{1-\beta}\int_{\partial K}h_L(\nu_K(x))|\nabla\varphi_K(x)|^2d\sigma=\frac{1+\beta}{1-\beta}\int_{\mathbb{S}^{N-1}}h_L(\xi)d\mu_K(\xi).$$}
\vskip 0.03in
Our aim of this section is to proving the following result.
\vskip 0.05in
{\em{\bf Theorem 4.1:} Let $K, L\in\mathcal{K}^N$ and $K_{tL}=K+tL$ for $t>0$. Then there holds
$$\frac{d}{dt}F(K_{tL}^\circ)_{t=0}=\frac{1+\beta}{1-\beta}\int_{\partial K}h_L(\nu_K(x))|\nabla\varphi_K(x)|^2d\sigma=\frac{1+\beta}{1-\beta}\int_{\mathbb{S}^{N-1}}h_L(\xi)d\mu_K(\xi).$$}
\vskip 0.03in
{\bf Proof:} At first, we prove Theorem 4.1 for $K\in\mathcal{K}_0^N$ and $L$ being of class $C^2$. To this end, we define $K_n=\{x\in K: \varphi_K(x)\geq\frac{1}{n}\}$. It is easy to see that $K_n$ is of class $C^2$ due to the smoothness of $\varphi_K$ in $K^\circ$ and $d_H(K_n,K)\to 0$ as $n\to\infty$. Therefore Theorem HSX can be used to conclude that for any $t>0$, we have
$$\frac{d}{dt}F((K_n+tL)^\circ)=\frac{1+\beta}{1-\beta}\int_{\mathbb{S}^{N-1}}h_L(\xi)d\mu_{K_n+tL}(\xi).$$
This implies that
\begin{equation}\label{eq41}
F((K_n+tL)^\circ)-F(K_n^\circ)=\frac{1+\beta}{1-\beta}\int_0^{t}\int_{\mathbb{S}^{N-1}}h_L(\xi)d\mu_{K_n+\tau L}(\xi)d\tau.
\end{equation}
Let $f_n(\tau)=\int_{\mathbb{S}^{N-1}}h_L(\xi)d\mu_{K_n+\tau L}(\xi)$. Then, by Theorem 3.1, we have
$$\lim\limits_{n\to\infty}f_n(\tau)=\int_{\mathbb{S}^{N-1}}h_L(\xi)d\mu_{K+\tau L}(\xi)$$
due to $d_H(K_n+\tau L,K+\tau L)\to 0$ as $n\to\infty$. Moreover, by Lemma 2.1, we can see that there exists a constant $C$ depending only on $K$ and $L$ such that $|f_n(\tau)|\leq C$. 
Therefore, it follows from the dominated convergence theorem that
\begin{equation}\label{eq42}
\lim\limits_{n\to\infty}\frac{1+\beta}{1-\beta}\int_0^{t}\int_{\mathbb{S}^{N-1}}h_L(\xi)d\mu_{K_n+\tau L}(\xi)d\tau=\frac{1+\beta}{1-\beta}\int_0^t\int_{\mathbb{S}^{N-1}}h_L(\xi)d\mu_{K+\tau L}(\xi)d\tau.
\end{equation}
By the continuity of $F$ (see Proposition 2.6), we have
\begin{equation}\label{eq43}
\lim\limits_{n\to\infty}(F((K_n+tL)^\circ)-F(K_n^\circ))=F((K+tL)^\circ)-F(K^\circ).
\end{equation}
By sending $n$ to infinity in equation (\ref{eq41}) and taking (\ref{eq42}), (\ref{eq43}) into account, we get
$$F((K+tL)^\circ)-F(K^\circ)=\frac{1+\beta}{1-\beta}\int_0^t\int_{\mathbb{S}^{N-1}}h_L(\xi)d\mu_{K+\tau L}(\xi)d\tau.$$
By the weak continuity of $\mu_\Omega$, we can see that $f(\tau)=\int_{\mathbb{S}^{N-1}}h_L(\xi)d\mu_{K+\tau L}(\xi)d\tau$ is continuous. Therefore, $F(K+tL)$ is differentiable and
\begin{equation}\label{eq44}
\frac{d}{dt}F(K_{tL}^\circ)=\frac{1+\beta}{1-\beta}\int_{\mathbb{S}^{N-1}}h_L(\xi)d\mu_{K+tL}(\xi)d\tau.
\end{equation}
For general $K, L\in\mathcal{K}$, we choose a sequence of $C^2$ convex bodies $L_n$ such that $d_H(L_n,L)\to 0$ as $n\to\infty$. By (\ref{eq44}), we have
\begin{equation}\label{eq45}
F((K+tL_n)^\circ)-F(K^\circ)=\frac{1+\beta}{1-\beta}\int_0^t\int_{\mathbb{S}^{N-1}}h_{L_n}(\xi)d\mu_{K+\tau L_n}(\xi)d\tau.
\end{equation}
By the uniform convergence of $h_{L_n}$ to $h_L$ and the weak continuity of the measure, we can conclude that
\begin{equation}\label{eq46}
\lim\limits_{n\to\infty}\int_0^t\int_{\mathbb{S}^{N-1}}h_{L_n}(\xi)d\mu_{K+\tau L_n}(\xi)d\tau=\int_0^t\int_{\mathbb{S}^{N-1}}h_L(\xi)d\mu_{K+\tau L}(\xi)d\tau.
\end{equation}
In fact, by the uniform convergence of $h_{L_n}$ to $h_L$ and the finiteness of $\mu_{K+\tau L_n}$, we have
$$\lim\limits_{n\to\infty}\int_0^t\int_{\mathbb{S}^{N-1}}(h_{L_n}(\xi)-h_L)d\mu_{K+\tau L_n}(\xi)d\tau=0.$$
By the weak continuity of the measure $\mu_{K^\circ}$, we have
$$\lim\limits_{n\to\infty}\int_{\mathbb{S}^{N-1}}h_L(\xi)d\mu_{K+\tau L_n}(\xi)d\tau=\int_0^t\int_{\mathbb{S}^{N-1}}h_L(\xi)d\mu_{K+\tau L}(\xi)d\tau.$$
Therefore
\begin{equation*}
\begin{array}{ll}
\lim\limits_{n\to\infty}\int_0^t\int_{\mathbb{S}^{N-1}}h_{L_n}(\xi)d\mu_{K+\tau L_n}(\xi)d\tau&=\lim\limits_{n\to\infty}\int_0^t\int_{\mathbb{S}^{N-1}}(h_{L_n}(\xi)-h_L)d\mu_{K+\tau L_n}(\xi)d\tau\\
&+\lim\limits_{n\to\infty}\int_0^t\int_{\mathbb{S}^{N-1}}h_L(\xi)d\mu_{K+\tau L_n}(\xi)d\tau\\
&=\int_0^t\int_{\mathbb{S}^{N-1}}h_L(\xi)d\mu_{K+\tau L}(\xi)d\tau.
\end{array}
\end{equation*}
From this and (\ref{eq45}), we can deduce that
$$F((K+tL)^\circ)-F(K^\circ)=\frac{1+\beta}{1-\beta}\int_0^t\int_{\mathbb{S}^{N-1}}h_L(\xi)d\mu_{K+\tau L}(\xi)d\tau.$$
This leads to our final conclusion
$$\frac{d}{dt}F(K_{tL}^\circ)|_{t=0}=\frac{1+\beta}{1-\beta}\int_{\mathbb{S}^{N-1}}h_L(\xi)d\mu_{K}(\xi).$$

\vskip 0.05in
\section*{5. Solution to Minkowski problem and its applications}

\setcounter{section}{5}

\setcounter{equation}{0}

\noindent

Let $K\in\mathcal{K}^N$, $\varphi_K(x)$ be the unique solution to problem (\ref{eq221}) in $K^\circ$ and $\mu_{K^\circ} $ be the Borel measure given by Definition 2.3. In this section, we try to solving Minkowski type problem associated with $\mu_{K^\circ}$. More precisely, for a given positive finite Borel measure $\mu$, we try to find a $K\in\mathcal{K}^N$ such that $\mu_{K^\circ}=\mu$. By the translation invariant of problem (\ref{eq221}), we can easily see that
$$\int_{\mathbb{S}^{N-1}}\xi d\mu_{K^\circ}(\xi)=0.$$
Therefore, a necessary condition for the solvability of our Minkowski type problem is that the given measure $\mu$ satisfies
$$\int_{\mathbb{S}^{N-1}}\xi d\mu(\xi)=0.$$
Also, by the translation invariant of problem (\ref{eq221}), we can consider our problem only in the class of $\mathcal{K}_0^N$. For any $K\in\mathcal{K}_0^N$, we define the mean width $W_\mu(K)$ of $K$ associated with a positive measure $\mu$ on $\mathbb{S}^{N-1}$ by
$$W_\mu(K)=\frac{2}{\mu(\mathbb{S}^{N-1})}\int_{\mathbb{S}^{N-1}}h_K(\xi)d\mu$$
Since the measure $\mu_{K^\circ}$ can be realized by the variation of domain functional $F(K^\circ)=\int_{K^\circ}|\nabla\varphi_K|^2dx$, the solution of our Minkowski problem can be realized by a stationary point of $F(K^\circ)$. This leads us to consider the following maximization problem
\begin{equation}\label{eq51}
m=\sup\limits_{K\in\mathcal{K}^N_0}\{F(K^\circ): W_\mu(K)=1\}.
\end{equation}
Let $B_{r}(0)=\{x: |x|<r\}$. Then $W_\mu(\overline{B}_{1/2}(0))=1$. Therefore, we have $m\geq F(B_{1/2}(0))>0$. The main result of this section is the following result.
\vskip 0.03in
{\em{\bf Theorem 5.1:} Let $\mu$ be a positive finite Borel measure on $\mathbb{S}^{N-1}$ such that
$$\int_{\mathbb{S}^{N-1}}\xi d\mu(\xi)=0\ \ \mbox{and}\ \ \int_{\mathbb{S}^{N-1}}|\theta\cdot\xi|d\mu(\xi)>0\ \ \mbox{for all}\ \ \theta\in\mathbb{S}^{N-1}.$$
Then, the $m$ defined by (\ref{eq51}) can be achieved by a convex body $K_m\in\mathcal{K}_0^N$.}
\vskip 0.03in
{\bf Proof:} Let $\{K_n\}_{n=1}^\infty$ be a maximizer sequence of (\ref{eq51}) in $\mathcal{K}_0^N$. That is $K_n\in\mathcal{K}_0^N$ verifies
\begin{equation}\label{eq52}
W_\mu(K_n)=\frac{2}{\mu(\mathbb{S}^{N-1})}\int_{\mathbb{S}^{N-1}}h_{K_n}(\xi)d\mu=1.
\end{equation}
\begin{equation}\label{eq53}
\lim\limits_{n\to\infty}F(K_n^\circ)=m.
\end{equation}
We claim that there are two positive number $R_{in}$ and $R_{out}$ such that
\begin{equation}\label{eq54}
B_{R_{in}}(0)\subset K_n\subset B_{R_{out}}(0).
\end{equation}
In fact, by the assumptions of Theorem 5.1, we know that the classical Minkowski problem associated with $\mu$ is uniquely solvable. That is, there is a convex body $K_{C}^\mu\in\mathcal{K}_0^N$ such that $d\sigma_{K_{C}^\mu}=g_*(d\sigma)=d\mu$ with $g: \partial K_{C}^\mu\to\mathbb{S}^{N-1}$ being the Gauss mapping of $K_{C}^\mu$ and $d\sigma$ being the surface element of $\partial K_{C}^\mu$. Since $K_C^\mu$ is convex body containing the origin, there exists a ball $B_r(0)$ with radius $r>0$ and centered at the origin such that $B_r(0)\subset K_C^\mu$. This implies that there is a positive constant $C(r)$ depending only on $r$ such that
$$\int_{\mathbb{S}^{N-1}}|\theta\cdot\xi|d\mu_{K_C^\mu}(\xi)>C(r)>0\ \ \mbox{for all}\ \ \theta\in\mathbb{S}^{N-1}.$$
Noting $d\sigma_{K_{C}^\mu}=d\mu$, we get
 \begin{equation}\label{eq55}
 \int_{\mathbb{S}^{N-1}}|\theta\cdot\xi|d\mu(\xi)>C(r)>0\ \ \mbox{for all}\ \ \theta\in\mathbb{S}^{N-1}.
 \end{equation}
For arbitrary $n$, let $D_n$ denote the diameter of $K_n$. No loss of generality, we may assume that for some $\theta\in\mathbb{S}^{N-1}$ there holds $\pm\frac{D_n}{2}\theta\in\partial K_n$. Then, for every $\xi\in\mathbb{S}^{N-1}$, we have
$$h_{K_n}(\xi)=\sup\limits_{x\in\partial K_n}x\cdot\xi\geq\max\{\pm\frac{D_n}{2}\theta\cdot\xi\}=\frac{D_n}{2}\theta\cdot\xi.$$
Combining this with (\ref{eq52}) and (\ref{eq55}) together, we obtain
$$\mu(\mathbb{S}^{N-1})=2\int_{\mathbb{S}^{N-1}}h_{K_n}(\xi)d\mu\geq D_n\int_{\mathbb{S}^{N-1}}\theta\cdot\xi d\mu(\xi)\geq D_nC(r).$$
This leads to $D_n\leq\frac{\mu(\mathbb{S}^{N-1})}{C(r)}$ for arbitrary $n$. Therefore, we can take $R_{out}=\frac{\mu(\mathbb{S}^{N-1})}{C(r)}$.

To estimate $R_{in}$, we recall a concept of John ellipsoid of a convex body $K$, A John ellipsoid of a convex body $K$ is the ellipsoid of minimal volume containing $K$. John's famous theorem says that for any convex body $K$ there is a John ellipsoid $\mathcal{E}_K$ such that
$$\frac{1}{N}\mathcal{E}_K\subset K\subset\mathcal{E}_K$$
provided that the center of $\mathcal{E}_K$ is properly chosen. Let $\mathcal{E}_n$ be the John ellipsoid of $K_n$ and $0<a_n^1\leq a_n^2\leq\cdots\leq a_n^N$ the lengths of the semi-major axes of $\mathcal{E}_n$. Then there holds
$$a_n^N\leq NR_{out}\ \ \mbox{and}\ \ V(K_n)\leq\omega_{N}\prod\limits_{j=1}^Na_n^j\leq a_n^1\omega_NN^{N-1}R_{out}^{N-1}.$$ 
with $\omega_N$ being the volume of unit ball in $\mathbb{R}^N$. By Lemma 2.1, there is a positive constant $C_1$ independent of $n$ such that $|\nabla\varphi_{K_n}(x)|\leq C_1$.
Therefore
$$F(K_n^\circ)=\int_{K_n^\circ}|\nabla\varphi_{K_n}(x)|^2dx\leq C_1V(K_n)\leq C_1a_n^1\omega_NN^{N-1}R_{out}^{N-1}.$$
Noting $\lim\limits_{n\to\infty}F(K_n^\circ)=m\geq F(B_{1/2}(0))>0$, we conclude that there is positive number $\rho$ independent of $n$ such that $a_n^1\geq\rho>0$. This implies that there exists a positive number $R_{in}$ such that $B_{R_{in}}(0)\subset K_n$ So the claim (\ref{eq54}) is concluded. 

From claim (\ref{eq54}) and the Blaschke's selection theorem, we may assume, up to a subsequence , that $K_n$ converges to $K\in\mathcal{K}_0^N$ in the Hausdorff metric. By Proposition 2.6 and the uniform convergence of $h_{K_n}$ to $h_{K}$, we have
$$W_\mu(K)=1\quad F(K^\circ)=m.$$
This completes the proof of Theorem 5.1.
\vskip 0.03in
As a corollary of Theorem 5.1, we have the following result.
\vskip 0.05in
{\em{\bf Theorem 5.2} Let $\mu$ be a positive finite Borel measure on $\mathbb{S}^{N-1}$ such that
$$\int_{\mathbb{S}^{N-1}}|\theta\cdot\xi|d\mu(\xi)>0\ \ \mbox{for all}\ \ \theta\in\mathbb{S}^{N-1}.$$
Then, there is a convex body $K\in\mathcal{K}^N$ such that $\mu_K=\mu$ if and only if $\int_{\mathbb{S}^{N-1}}\xi d\mu(\xi)=0$.}}
\vskip 0.03in
{\bf Proof:} The necessary part follows immediately from the translation invariant of problem (\ref{eq221}). What we should do is to proving the sufficient part. Let $K$ be the maximizer of maximization problem (\ref{eq51}). For any $L\in\mathcal{K}_0^N$, let $K_t=K+tL$. Since support function is linear with respect to Minkowski addition, we have
$$\frac{d}{dt}W_\mu(K_t)|_{t=0}=\frac{2}{\mu(\mathbb{S}^{N-1})}\int_{\mathbb{S}^{N-1}}h_L(\xi)d\mu(\xi).$$
By Theorem 4.1, we have
$$\frac{d}{dt}F(K_t^\circ)|_{t=0}=\frac{1+\beta}{1-\beta}\int_{\mathbb{S}^{N-1}}h_L(\xi)d\mu_{K}(\xi).$$
It follows from the Lagrange multiplier method that there exists a positive number $\lambda$ such that
\begin{equation}\label{eq56}
\int_{\mathbb{S}^{N-1}}h_L(\xi)d\mu_{K}(\xi)=\lambda\int_{\mathbb{S}^{N-1}}h_L(\xi)d\mu(\xi).
\end{equation}
To conclude $\mu_K=\lambda\mu$, we must prove
\begin{equation}\label{eq57}
\begin{array}{ll}
\int_{\mathbb{S}^{N-1}}f(\xi)d\mu_{K}(\xi)=\lambda\int_{\mathbb{S}^{N-1}}f(\xi)d\mu(\xi)&\forall f\in C(\mathbb{S}^{N-1}).
\end{array}
\end{equation}
For $f\in C^2(\mathbb{S}^{N-1})$, it follows from Lemma 1.7.8 of \cite{S2014} that there are two convex bodies $L_1, L_2\in\mathcal{K}_0^N$ such that
\begin{equation}\label{eq58}
f=h_{L_1}-h_{L_2}
\end{equation}
By making use of (\ref{eq56}), we have
\begin{equation}\label{eq59}
\begin{array}{ll}
\int_{\mathbb{S}^{N-1}}f(\xi)d\mu_{K}(\xi)&=\int_{\mathbb{S}^{N-1}}h_{L_1}(\xi)d\mu_{K}(\xi)-\int_{\mathbb{S}^{N-1}}h_{L_2}(\xi)d\mu_{K}(\xi)\\
&=\lambda\int_{\mathbb{S}^{N-1}}h_{L_1}(\xi)d\mu(\xi)-\lambda\int_{\mathbb{S}^{N-1}}h_{L_2}(\xi)d\mu(\xi)\\
&=\lambda\int_{\mathbb{S}^{N-1}}f(\xi)d\mu(\xi).
\end{array}
\end{equation}
For any given $f\in C(\mathbb{S}^{N-1})$, we can choose a sequence of functions $f_n\in C^2(\mathbb{S}^{N-1})$ such that
$$\lim\limits_{n\to\infty}\|f_n-f\|_{L^\infty(\mathbb{S}^{N-1})}=0.$$
By making use of (\ref{eq59}), we obtain
\begin{equation}\label{eq510}
\int_{\mathbb{S}^{N-1}}f_n(\xi)d\mu_{K}(\xi)=\lambda\int_{\mathbb{S}^{N-1}}f_n(\xi)d\mu(\xi).
\end{equation}
Because $\mu_K$ and $\mu$ are finite, we can easily see that
$$\lim\limits_{n\to\infty}\int_{\mathbb{S}^{N-1}}f_n(\xi)d\mu_K=\int_{\mathbb{S}^{N-1}}f(\xi)d\mu_K,\quad\lim\limits_{n\to\infty}\int_{\mathbb{S}^{N-1}}f_n(\xi)d\mu=\int_{\mathbb{S}^{N-1}}f(\xi)d\mu.$$
Taking limit as $n\to\infty$ in equation (\ref{eq510}), we have
\begin{equation}\label{eq57}
\begin{array}{ll}
\int_{\mathbb{S}^{N-1}}f(\xi)d\mu_{K}(\xi)=\lambda\int_{\mathbb{S}^{N-1}}f(\xi)d\mu(\xi)&\forall f\in C(\mathbb{S}^{N-1}).
\end{array}
\end{equation}
Therefore $\mu_K=\lambda\mu$. Since $\mu_K$ is $\alpha=N+\frac{3\beta+1}{1-\beta}$ homogeneous, we see that $\widetilde{K}=\lambda^{-\alpha} K$ satisfies $\mu_{\widetilde{K}}=\mu$. This completes the proof of Theorem 5.2.
\vskip 0.05in
{\em{\bf Remark 5.3:} By Theorem 11 of \cite{C2005}, the solution of $\mu_K=\mu$ is unique up to a translation.}
\vskip 0.05in
Let $H^{N-1}$ be the $(N-1)$-dimensional Hausdorff measure. To prove Theorem 1.3, we first recall a result of \cite{Fr} due to I. Fragal$\grave{a}$ which can be stated as
\vskip 0.03in
{\em{\bf Lemma 5.4:} Let $F: \mathcal{K}_0^N\to\mathbb{R}_+$ be a  Brunn-Minkowski functional, whose definition is given in \cite{Fr}, of degree $\alpha$. If $K$ is a stationary point of the functional
$$\Psi(K)=\frac{F^{\frac{1}{\alpha}}(K)}{M_{H^{N-1}}(K)}$$
in the sense that $\frac{d}{dt}\Psi((1-t)K+tL)|_{t=0}=0$ for any $L\in\mathcal{K}_0^N$, then $K$ is a ball.}

It is easy to see that the functional $F(K)$ which we use in Theorem 5.1 is a Brunn-Minkowski functional of degree $\alpha=N+\frac{2(1+\beta)}{1-\beta}$. Let $\mu=H^{N-1}$. By Theorem 5.2, the equation $\mu_K=H^{N-1}$ has a unique solution $\overline{K}$. By Theorem 5.1, we can verify that $\overline{K}$ is a stationary point of $\Psi(K)=\frac{F^{\frac{1}{\alpha}}(K)}{M_{H^{N-1}}(K)}$ in the sense of Lemma 5.4. Therefore, Lemma 5.4 ensures that $\overline{K}$ is a ball. This combine with the uniqueness result leads to the following isoperimetric result.
\vskip 0.05in
{\em{\bf Theorem 5.5:} Let $H^{N-1}$ denote the $(N-1)$-dimensional Hausdorff measure and $B$ denote a ball in $\mathbb{R}^N$. Then, for any bounded convex domain $K$, there holds $F(K)\leq F(B)$ provided that $W_{H^{N-1}}(K)=W_{H^{N-1}}(B)$. Moreover, $F(K)=F(B)$ and $W_{H^{N-1}}(K)=W_{H^{N-1}}(B)$ imply $\Omega=B$ up to a translation.} 
\vskip 0.05in
Since $W_{H^{N-1}}$ is homogeneous of order one, we see that $W_{H^{N-1}}(\frac{W_{H^{N-1}}(B_1)}{W_{H^{N-1}}(K)}K)=W_{H^{N-1}}(B_1)=2$. 
Therefore, for any bounded convex domain $K$
$$F(K)\leq\frac{F(B_1)}{2^\alpha}W_{H^{N_1}}^\alpha(K)$$
with $\alpha=N+\frac{2(1+\beta)}{1-\beta}$, and the equality holds if and only if $K$, up to a translation, is a ball. 
\vskip 0.03in
{\em{\bf Remark 5.6:} When we post our work on ArXive.org, Frank tell us that a result similar to Theorem 5.5 for general Lipschitz domains is obtain in their recent work \cite{FL}.}

\end{document}